\documentclass[11pt]{article}
\usepackage{amsmath,amssymb}
\usepackage{pstricks}
\usepackage{pst-node}
\usepackage{pst-poly}
\usepackage{multido}
\usepackage{tikz}
\newtheorem{result}{Theorem}

\newtheorem{support}{Lemma}
\newtheorem{propo}{Proposition}

\newcommand{\qed}{%
\ifmmode % if math mode, assume display: omit penalty etc.
\else \leavevmode\unskip\penalty9999 \hbox{}\nobreak\hfill \fi
\quad\hbox{\qedsymbol}}
\newcommand{\openbox}{\leavevmode \hbox to.77778em{%
\hfil\vrule
\vbox to.675em{\hrule width.6em\vfil\hrule}%
\vrule\hfil}}
\newcommand{\qedsymbol}{\openbox}
\newcommand{\showgrid}{}
\newcommand{\gridon}{\renewcommand{\showgrid}{\psset{subgriddiv=1,griddots=10,gridlabels=6pt}\psgrid}}
\gridon
\begin{document}

\begin{center} {\LARGE  Packing 1-Plane Hamiltonian Cycles in Complete Geometric Graphs} \end{center}
\vskip8pt

\centerline{\small Hazim Michman Trao, Adem Kilicman \ and \ Niran Abbas Ali}

\begin{center}
\itshape\small  Department of Mathematics, \\ Universiti Putra Malaysia, 43400 Serdang, Malaysia,  \\

\end{center}

\begin{abstract} Counting the number of Hamiltonian cycles that are contained in a geometric graph is {\bf \#P}-complete even if the graph is known to be planar \cite{lot:refer}. A relaxation for problems in plane geometric graphs is to allow the geometric graphs to be 1-plane, that is, each of its edges is crossed at most once. We consider the following question: For any set $P\/$ of $n\/$ points in the plane, how many 1-plane Hamiltonian cycles can be packed into a complete geometric graph $K_n\/$? We investigate the problem by taking two different situations of $P\/$, namely, when $P\/$ is in convex position, wheel configurations position. For points in general position we prove the lower bound of $k-1\/$ where $n=2^{k}+h\/$ and $0\leq h <2^{k}\/$. In all of the situations, we investigate the constructions of the graphs obtained.
\end{abstract}

\vspace{1mm}
\section{Introduction}

Let $P\/$ be a set of $n\/$ points in general position in the plane with no three points being collinear. A {\em geometric graph\/} is a graph $G = (P,E)$ that consists of a set of vertices $P$, which are points in the plane, and a set of edges, $E$, which are straight-line segments whose endpoints belong to $P$.  A {\em complete\/} geometric graph $K_n\/$ is a geometric graph on a set $P\/$ of $n\/$ points that has an edge joining every pair of points in $P\/$. Two edges are disjoint if they have no point in common. Two subgraphs are {\em edge-disjoint\/} if they do not share any edge.\

A geometric graph is said to be \textit{plane} (or non-crossing) if its edges do not cross each other. A geometric graph is said to be \textit{1-plane} if every edge is allowed to have at most one crossing. Note that the terms plane graph and 1-plane graph refer to a geometric object, while to be planar or 1-planar are properties of the underlying abstract graph.\

By an \emph{edge packing} of a graph $G\/$ we mean a set of edge-disjoint subgraphs of $G\/$. By an \emph{edge partition} of $G\/$ we mean an edge packing of $G\/$ with no edge left over, that is the union of all subgraphs in the packing is equal to $G\/$. Dor and Tarsi \cite{dt:refer} proved that the problem of partitioning a given graph $G\/$ is NP-complete.\

It is often useful to restrict the subgraphs of $G\/$ to a certain class or property. Among all subgraphs of $K_n\/$, plane spanning trees, plane Hamiltonian cycles or paths, and plane perfect matchings,  are of interest \cite{aam:refer, abd:refer, aic:refer, bhrw:refer, nm:refer} i.e., one may look for the maximum number of these subgraphs that can be packed into $K_n\/$. For instance, a long-standing open question is to determine if the edges of $K_n\/$, where $n\/$ is even, can be partitioned into $\frac{n}{2}$ plane spanning trees? Bernhart and Kanien \cite{bk:refer} give an affirmative answer for the problem when the points in convex position. Bose et al. \cite{bhrw:refer} proved that every complete geometric graph $K_n\/$ can be partitioned into at most $n-\sqrt{\frac{n}{12}}\/$ plane trees. Aichholzer et al. \cite{aic:refer} showed that $\Omega(\sqrt{n})$ plane spanning trees can be packed into $K_n\/$. Recently, Biniaz et al. \cite{bbms:refer} showed that at least $\lceil \log_2 n\rceil -1\/$ plane perfect matchings can be packed into $K_n\/$. \

A \textit{Hamiltonian cycle} is a cycle in a graph that passes through every vertex exactly once, except for the vertex that is both the beginning and end, which is visited twice. Finding a Hamiltonian cycle in a graph is {\bf NP}-complete even if the graph is known to be planar \cite{gjt:refer}. Moreover, counting the number of Hamiltonian cycles that are contained in a graph is {\bf \#P}-complete even if the graph is known to be planar \cite{lot:refer}. \

 The problem of counting the number of plane Hamiltonian cycles (not necessary edge-disjoint) on a given graph considered in \cite{aki:refer, dsst:refer, ssw:refer, sw:refer} and many others. In particular, Newborn and Moser \cite{nm:refer} asked for the maximal number of plane Hamiltonian cycles; authors give an upper bound of $2\cdot6^{n-2}\lfloor\frac{n}{2}\rfloor!$, and conjecture that it should be of the form $c^n$, for some constant $c$. Motzkin \cite{m:refer} proved that a set of $n\/$ points in the plane has at most $O(86:81^n)$ plane Hamiltonian cycles.\

A relaxation for problems in plane geometric graphs is to allow the geometric graphs to be 1-plane.  The problem of finding 1-plane Hamiltonian alternating cycle or path studied in many papers (see  \cite{kky:refer}, \cite{cggst:refer}). Claverol et al. \cite{cgh:refer} studied the 1-plane character. They showed that one can always obtain a 1-plane Hamiltonian alternating cycle on a point set in convex position and on a double chain.\

Since finding more than one edge-disjoint plane Hamiltonian cycle for a given set of points is not always possible to achieve, We relax the constraint on the Hamiltonian cycles from being plane to being 1-plane and study the following problem: For any set of $n\/$ points in the plane, how many 1-plane Hamiltonian cycles can be packed into a complete geometric graph $K_n\/$?\

For simplicity, we write 1-PHC to refer to a 1-plane Hamiltonian cycle.\

\subsection{Results}

We study the problem of packing 1-PHCs into the complete geometric graph $K_n\/$ for a given set of $n\/$ points in the plane. Since the complete graph $K_n\/$ on $n\/$ vertices has $n(n-1)/2\/$ edges and a Hamiltonian cycle has $n\/$ edges, therefore, the number of edge-disjoint Hamiltonian cycles in $K_n\/$ cannot exceed $(n-1)/2\/$.\

In Section 2, we show that $\lfloor\frac{n}{3}\rfloor\/$ is a tight bound for the number of 1-PHCs that can be packed into $K_n\/$ for any given set in convex position. In Section 3, we show that for a set of points in regular wheel configuration, $\lfloor\frac{n-1}{3}\rfloor\/$ edge-disjoint 1-PHCs can be packed into $K_n\/$ and this bound is tight. In the latter portion of this paper, point sets in general position are considered. We know that for $n\geq 3\/$, and by a minimum weight Hamiltonian cycle in $K_n\/$, a trivial lower bound of $1\/$ is obtained since it is a plane cycle. In Section 4, we present an algorithm (Algorithm $A\/$) to draw a 1-PHC for any set of points in general position in the plane. As the main result in this paper, we prove that there are at least $k-1\/$ 1-PHCs that can be packed into $K_n\/$, where $n=2^{k}+h\/$ and $0\leq h <2^{k}\/$.

\vspace {1mm} Throughout this paper, for simplicity, we consider all vertices in counter-clockwise order.

\section{1-PHCs for point sets in convex position}
In this section, we study the problem of packing 1-PHCs on a well-known restricted position of a point set which is the convex position. We will show that for any point set $P\/$ in convex position, there are at most $\lfloor \frac{n}{3}\rfloor$ edge-disjoint 1-PHCs that can be packed into $K_n\/$ and this bound is tight (Theorem \ref{the1}).

\vspace {1mm} Suppose $P= \{v_0, v_1, v_2, \ldots, v_{n-1}\}\/$ is a set of $n\/$ points in convex position. Let $G\/$ be a geometric graph on $P\/$. Edges of the form $v_iv_{i+1}\/$, $i=0,1,2,  \ldots, n-1\/$, are called the {\em boundary edges\/} in $G\/$. A non-boundary edge in $G\/$ is called a {\em diagonal edge\/}. \

\begin{propo} \label{pro2}
Let $P\/$ be a set of $n\/$ points in convex position in the plane where $n\geq 3\/$. Suppose $C\/$ is a 1-PHC on $P\/$ that has a diagonal edge $v_iv_j\/$ that divides $P\/$ into two parts, both including $v_i\/$ and $v_j\/$. Then the following statements hold:\

(1) If any part has an odd number of vertices, then $C\/$ has at least one boundary edge on this part.

(2) If any part has an even number of vertices, then $C\/$ has at least two boundary edges on this part.

\end{propo}

\vspace{1mm}  \noindent
{\bf Proof:} Let $C\/$ be a 1-PHC on $P\/$. Assume that $C\/$ contains the diagonal edge $v_iv_j\/$ that divides $P\/$ into two parts $P_1\/$ and $P_2\/$.\

\vspace{1mm} Assume that $P_1=\{v_j,v_{j+1},...,v_i\}\/$ and $|P_1|\/$ is odd. By induction on $|P_1|\/$, if $|P_1|=3\/$, then $P_1=\{v_j,v_{j+1},v_i\}\/$ and either the boundary edge $v_jv_{j+1}\/$ or $v_iv_{j+1}\/$ in $C\/$; otherwise, there are two crossings, a contradiction.\

\vspace{1mm} Assume that $|P_1|\geq 5\/$ is odd and the proposition is true when $m < |P_1|$, and $m$ is an odd number.\

\vspace{1mm} We claim that either $v_{i+1}v_{j+1}\/$ or $v_{i-1}v_{j-1}\/$ is an edge in $C\/$. To prove our claim, Suppose that neither $v_{i+1}v_{j+1}\/$ nor $v_{i-1}v_{j-1}\/$ are in $C\/$.\

\vspace{1mm} Since $C\/$ is a 1-PHC, then $C$ contains an edge $v_lv_k \notin \{v_{i+1}v_{j+1}, v_{i-1}v_{j-1}\}\/$ crosses $v_iv_j$ where $k\in P_1-\{v_i, v_j\}$ and $l\in P_2-\{v_i, v_j\}$.\

\vspace{1mm} When $k \notin \{j+1, i-1\}$, then there is $p\in\{j+2, j+3, ..., k-1\}$ or $p\in\{k+1, k+2, ..., j-2\}$ such that the edge $v_pv_q$, that incident to $v_p$ and belongs to $C$, crosses either $v_iv_j$ or $v_lv_k$, a contradiction.

\vspace{1mm} When $k \in \{j+1, i-1\}$, without loss of generality assume that $k=j+1$ and then $l\neq i+1$ (by assumption). Thus, there is $p\in\{j+2, j+3, ..., i-1\}$ or $p\in\{i+1, i+2, ..., l-1\}$ such that the edge $v_pv_q$, that belongs to $C$, crosses either $v_iv_j$ or $v_lv_k$, a contradiction.\

\vspace{1mm} Therefore, either $v_{i+1}v_{j+1}\in E(C)\/$ or $v_{i-1}v_{j-1}\in E(C)\/$. Without loss of generality, assume that $v_{i+1}v_{j+1}\in E(C)\/$.\

\vspace{1mm} Let $C'\/$ be a subgraph of $C\/$ induced by $P_1\/$. Since  $v_{i+1}v_{j+1}\/$ not in $C'\/$, then $d(v_j)=d(v_{j+1})=1\/$ in $C'\/$.\

\vspace{1mm} Let $C_1=C'\cup \{v_jv_{j+1}\}\/$. It is clear that $C_1\/$ is a 1-PHC on $P_1\/$ since $v_jv_{j+1}\/$ is a boundary edge. Recall that $v_iv_j\/$ and $v_jv_{j+1}\/$ are boundary edges in $C_1\/$ but are not boundary edges in $C\/$.\

\vspace{1mm} When the boundary edge $v_iv_{i-1}\in E(C_1)\/$., the claim in the proposition is hold.\

\vspace{1mm} Hence, assume that the boundary edge $v_iv_{i-1}\notin E(C_1)\/$. Then there is a diagonal edge $v_iv_k\in E(C_1)\/$ divides $P_1\/$ into two parts $P_{1,1}=\{v_i, v_{i-1}, ..., v_k\}\/$ and $P_{1,2}=\{v_i, v_{j}, v_{j+1}, ..., v_k\}\/$. \

\vspace{1mm} If $k=j+1\/$, there is a contradiction since $C_1\/$ is not a union of disjoint cycles.\

\vspace{1mm} Hence, either $k\in \{j+2,j+4,...,i-2\}\/$, then $P_{1,1}=\{v_k, v_{k+1},...,v_{i-1}, v_i\}\/$. Then clearly, $|P_{1,1}|\/$ is odd. By the induction hypothesis, $C_1\/$ has at least one boundary edge on $P_{1,1}\/$. Thus, $C\/$ has at least one boundary edge on $P_1\/$.\

\vspace{1mm} Or, $k\in \{j+3,j+5,...,i-3\}\/$, then $P_{1,2}=\{v_i,v_j,v_{j+1},...,v_{k}\}\/$. Then clearly, $|P_{1,1}|\/$ is odd.

\vspace{1mm} (*) Assume on the contrary that $C_1\/$ has the only two boundary edges $v_iv_j\/$ and $v_jv_{j+1}\/$ on $P_{1,2}\/$. Note that $\{v_{j+1},...,v_{k-1}\}\/$ has at least two vertices. However, $C_1\/$ matches all the vertices in $\{v_{j+1},...,v_{k-1}\}\/$ with at least two crossings with the edge $v_iv_k\/$ since $C_1\/$ has no boundary edge on $\{v_{j+1},...,v_{k-1}\}\/$; this is a contradiction (since $C_1\/$ is a 1-PHC). Thus, $C_1\/$ has at least one boundary edge different from $v_iv_j\/$ and $v_jv_{j+1}\/$ on $P_{1,2}\/$. This proves (1).\

\vspace{2mm} Assume that $|P_1|\/$ is even. By induction on $|P_1|\/$, if $|P_1|=4\/$, then $P_1=\{v_j,v_{j+1},v_{j+2},v_i\}\/$ such that either $v_jv_{j+1},v_{j+1}v_{j+2}\in E(C)\/$ or $v_iv_{j+2},v_{j+1}v_{j+2}\in E(C)\/$; otherwise, there is a contradiction since $C\/$ is a 1-PHC. Assume that $|P_1|\geq 6\/$ is even and the proposition is true when $m < |P_1|$, $m$ is even.\

\vspace{1mm} By repeating the same argument in (1), we conclude that either $v_{i+1}v_{j+1}\in E(C)\/$ or $v_{i-1}v_{j-1}\in E(C)\/$. Without loss of generality, assume that $v_{i+1}v_{j+1}\in E(C)\/$.\

\vspace{1mm} Let $C'\/$ be a subgraph of $C\/$ induced by $P_1\/$. It is clear that $d(v_j)=d(v_{j+1})=1\/$ in $C'\/$. Let $C_1=C'\cup \{v_jv_{j+1}\}\/$. Then $C_1\/$ is a 1-PHC on $P_1\/$ since $v_jv_{j+1}\/$ is a boundary edge. Recall that $v_iv_j\/$ and $v_jv_{j+1}\/$ are boundary edges in $C_1\/$ but are not boundary edges in $C\/$.\

\vspace{1mm} In the case that $C_1\/$ has the boundary edge $v_{i}v_{i-1}$, either $v_{i-1}v_{i-2}\in E(C)$ and then the claim in the proposition is true, or  $v_{i-1}v_{i-2}\notin E(C)$ and then the diagonal edge $v_{i-1}v_{k}\in E(C)$ for some $k\in\{j+2, j+3, ..., i-3\}$.

\vspace{1mm}Hence, the diagonal $v_{i-1}v_{k}$ divides $P_1\/$ into two parts $P_{1,1}=\{v_{k}, v_{k+1}, ..., v_{i-2}, v_{i-1}\}\/$ and $P_{1,2}=\{v_{i-1}, v_i, v_{j}, v_{j+1}, ..., v_k\}\/$. \

\vspace{1mm} $|P_{1,1}|\/$, whether odd or even, $C\/$ has at least one boundary edge on $P_{1,1}\/$ by part (1) or by induction, respectively.\

\vspace{1mm} In the case that $C_1\/$ does not have the boundary edge $v_{i}v_{i-1}$, then $v_{i}v_{k}\in E(C)$ for some $k\in\{j+2, j+3, ..., i-2\}$.

\vspace{1mm}Hence, the diagonal edge $v_{i}v_{k}$ divides $P_1\/$ into two parts $P_{1,1}=\{v_{k}, v_{k+1}, ..., v_{i-1}, v_{i}\}\/$ and $P_{1,2}=\{v_i, v_{j}, v_{j+1}, ..., v_k\}\/$. \

\vspace{1mm}  Now, either both $|P_{1,1}|\/$ and $|P_{1,2}|\/$ are odd, then $C\/$ has at least one boundary edge on $P_{1,1}\/$ by part (1) and it has at least one boundary edge different from $v_iv_j\/$ and $v_jv_{j+1}\/$ on $P_{1,2}\/$ by argument (*). \

\vspace{1mm}  Or, both $|P_{1,1}|\/$ and $|P_{1,2}|\/$ are even. By the induction hypothesis, $C_1\/$ has at least two boundary edges on $P_{1,1}\/$.\

\vspace{1mm} Thus, $C\/$ has at least two boundary edges on $P_1\/$. This completes the proof. \qed

\vspace{1mm} As a direct consequence of Proposition \ref{pro2}, we have the following lemma.

\begin{support} \label{lem1}

Let $P\/$ be a set of $n\/$ points in convex position in the plane where $n\geq 3\/$. Suppose $C\/$ is a 1-PHC on $P\/$. Then the following statements hold:\

(1) If $n$ is even, $C\/$ has at least two boundary edges.

(2) If $n$ is odd, $C\/$ has at least three boundary edges.

\end{support}

\vspace{1mm}  \noindent
{\bf Proof:} Let $C\/$ be a 1-PHC on a set $P\/$ of $n\/$ points. If all edges of $C\/$ are boundary edges, then the claim in the lemma holds. Thus, assume that $C\/$ contains a diagonal edge $v_iv_j\/$.\

\vspace{1mm} If $n\/$ is even, $v_iv_j\/$ divides $P\/$ into two parts, each part having an odd (even) number of vertices. Then by Proposition \ref{pro2}, $C\/$  has at least one boundary edge (two boundary edges) on each part.\

\vspace{1mm} If $n\/$ is odd, $v_iv_j\/$ divides $P\/$ into two parts, and one part has an odd number of vertices. By Proposition \ref{pro2}, $C\/$  has at least one boundary edge on this part, while $C\/$ has at least two boundary edges on the second part, which has an even number of vertices.\qed
\\

Suppose $G\/$ is a geometric graph on a set of convex position $P\/$ that has a diagonal edge $v_iv_j\/$. A boundary edge $v_kv_{k+1}\/$ is called {\em on the right side} of $v_iv_j\/$ if $i\leq k< j$ and {\em on the left side} of $v_iv_j\/$ if $j\leq k< i$. A diagonal edge is said to {\em have a boundary edge on each side} if there are two boundary edges, on left and right sides of the diagonal edge. A boundary edge $v_kv_{k+1}\/$ is called a {\em single boundary edge in $G\/$} if the two boundary edges $v_{k-1}v_k\/$ and $v_{k+1}v_{k+2}\/$ are not in $G\/$.\

\begin{propo} \label{pro3}
Let $P\/$ be a set of $n\/$ points in convex position in the plane where $n\geq 4\/$. Suppose $C\/$ is a 1-PHC on $P\/$. Then
the following statements hold:

(1) If $C\/$ has only two boundary edges $\{v_kv_{k+1}\}$ for $k\in \{r,s\}$, then $C\/$ has the edges $v_kv_{k+2}$ and $v_{k+1}v_{k-1}$.

(2) If $C\/$ has only three boundary edges, then $C\/$ has at least one single boundary edge $v_rv_{r+1}$ and the edges $v_rv_{r+2}$ and $v_{r+1}v_{r-1}$.

\end{propo}

\vspace{1mm}  \noindent
{\bf Proof:} Let $C$ be a 1-PHC on $P\/$ containing only two boundary edges $v_kv_{k+1}$ for $k\in \{r,s\}$. Assume on the contrary that at least one of the two edges $\{v_kv_{k+2}, v_{k+1}v_{k-1}\}\/$ is not in $C\/$ for some $k\in \{r,s\}$.\

\vspace{1mm} Without loss of generality, assume that $v_rv_{r+2}\notin E(C)\/$. Then $v_rv_i\/$ and $v_{r+1}v_j\in E(C)$ for some $r+3\leq i\leq r-2\/$ and $r+3\leq j\leq r-1\/$.\

\vspace{1mm} If $v_rv_{r+1}\/$ and $v_sv_{s+1}\/$ are in consecutive order, then all the remaining edges of $C\/$ are diagonal edges where $v_rv_{r+1}\/$ and $v_sv_{s+1}\/$ are on the same side of each of which. Hence, the other side of any diagonal edge does not contains any boundary edge, which is a contradiction with Proposition \ref{pro2}. Thus $v_kv_{k+1}$ is a single boundary edge for each $k\in \{r,s\}$ and then $i\neq r+1\/$ and $j\neq r+2\/$. \

\vspace{1mm}
If $v_rv_i\/$ and $v_{r+1}v_{j}$ are crossing, then there is a vertex $v_t\/$ where $r+2\leq t<i\/$ such that at least one of the two edges that incident to $v_t\/$ in $C\/$ crosses $v_rv_i\/$, which is a contradiction since $C\/$ is a 1-PHC.\

\vspace{1mm} If $v_rv_i\/$ and $v_{r+1}v_j\/$ are not crossing, then by Proposition \ref{pro2}, $C\/$ has at least one boundary edge on the left side of $v_rv_i\/$ and at least another boundary edge on the right side of $v_{r+1}v_j\/$ (both different from $v_rv_{r+1}\/$), which is a contradiction since $C\/$ has only two boundary edges. This proves (1).\

\vspace{1mm} Let $C\/$ contain three boundary edges $v_kv_{k+1}\/$ for $k\in \{r,s,t\}\/$. Suppose that no single boundary edge in $C\/$. That is, the boundary edges in $C\/$ are in consecutive order. But all the remaining edges of $C\/$ are diagonal edges where $v_kv_{k+1}\/$ for each $k\in \{r,s,t\}\/$ are on the same side of each of which. Hence, the other side of any diagonal edge does not contains any boundary edge, which is a contradiction with Proposition \ref{pro2}. Thus $C\/$ has at least one single boundary edge.\

\vspace{1mm} Assume on the contrary that, if $v_kv_{k+1}\/$ is a single boundary edge in $C\/$ for some $k\in \{r,s,t\}\/$, then either $v_kv_{k+2}\/$ or $v_{k+1}v_{k-1}\/$ or both are not in $C\/$.\

\vspace{1mm} Without loss of generality, assume that $v_rv_{r+1}\/$ is a single boundary edge in $C\/$ and $v_rv_{r+2}\notin E(C)$. Then $v_rv_i\/$ and $v_{r+1}v_j\/$ in $C\/$ where $r+3\leq i< r-1\/$ and $r+2< j\leq r-2\/$.\

\vspace{1mm} If $v_rv_i\/$ and $v_{r+1}v_j\/$ are crossing, then there is a vertex $v_t\/$ where $r+2\leq t<j$ such that at least one of the two edges that incident to $v_t\/$ in $C\/$ crosses $v_rv_i\/$, which is a contradiction since $C\/$ is a 1-PHC.\

\vspace{1mm} If $v_rv_i\/$ and $v_{r+1}v_j\/$ are not crossing, then by Proposition  \ref{pro2}, $C\/$ has at least one boundary edge on the set $\{v_{i+1},v_{i+2},...,v_r\}\/$ and at least one boundary edge on the set $\{v_{r+1},v_{r+2},...,v_j\}\/$; otherwise, there is a contradiction since $C\/$ has only three boundary edges.\

\vspace{2mm} This implies that there is a single boundary edge $v_sv_{s+1}\/$ where $i\leq s<r-1\/$ such that either $v_sv_{s+2}\/$ or $v_{s+1}v_{s-1}\/$ is not in $C\/$ (by assumption). Without loss of generality, assume that $v_sv_{s+2}\notin E(C)\/$. Let $v_sv_p\/$ and $v_{s+1}v_q\/$ in $C\/$.\

\vspace{1mm} Suppose that $p\in\{j,j+1,...,s-1\}\/$; otherwise, $v_sv_p\/$ crosses both $v_rv_i\/$ and $v_{r+1}v_j\/$, which is a contradiction since $C\/$ is a 1-PHC. Then by Proposition \ref{pro2}, $C\/$ has at least one boundary edge on a vertex set $\{v_j,v_{j+1},...,v_{s-1}\}\/$, which is a contradiction since $C\/$ has only three boundary edges. This completes the proof. \qed

\vspace{2mm}
Easily, one can see the Proposition \ref{pro3} still true when $n=3\/$, is odd with consecutive boundary edges.\

\vspace{1mm}
\begin{result}\label{the1}
Let $P\/$ be a set of $n\/$ points in convex position on the plane where $n\geq 3\/$. Then there exist $k\/$ edge-disjoint 1-PHCs $C_1, C_2, ... , C_k\/$ on $P\/$ that can be packed into $K_n\/$ where $k\leq \lfloor \frac{n}{3}\rfloor \/$.
\end{result}

\vspace{1mm}  \noindent
{\bf Proof:} (1) Let $n=2m\/$, which is even. Suppose $P=\{v_0, v_1, . . . , v_{n-1}\}\/$ is a set of $n\/$ points in convex position on the plane. By Lemma \ref{lem1}, every 1-PHC on $P\/$ contains at least two boundary edges. On the other hand, $P\/$ has $n\/$ boundary edges; that is, the number of 1-PHCs does not exceed $n/2\/$.\

\vspace{1mm} We claim that if $C_i\/$ and $C_j\/$ are two edge-disjoint 1-PHCs each having only two boundary edges, then any boundary edge of $C_i\/$ can not be in consecutive order with a boundary edge of $C_j\/$. To prove our claim, assume on the contrary that $v_rv_{r+1}\in E(C_i)\/$ and $v_{r+1}v_{r+2}\in E(C_j)\/$. By Proposition \ref{pro3}, the edge $v_rv_{r+2}\in E(C_i)\cap E(C_j)\/$, which is a contradiction since $C_i\/$ and $C_j\/$ are edge-disjoint 1-PHCs ($E(C_i)\cap E(C_j)=\phi\/$, where $i\neq j\/$).\

\vspace{1mm} Note that a single boundary edge in $C_i\/$ can be adjacent to any two consecutive boundary edges in $C_j\/$ where $i\neq j\/$; that is, $C_i \cup C_j\/$ can have three boundary edges in consecutive order. Therefore, the number of 1-PHCs that can be packed into $K_n$ is at most $\lfloor \frac{n}{3}\rfloor\/$ and this bound is tight.\

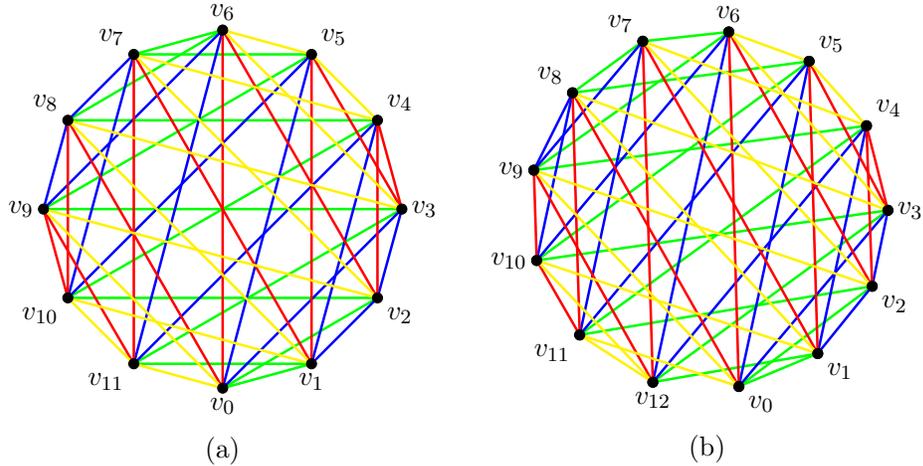
\begin{figure}[htb]
\centering
\resizebox{13cm}{!}{
\begin{minipage}{.45\textwidth}
\begin{tikzpicture}

           %points
        \coordinate (v0) at (0,-2.5);\filldraw[black] (v0) circle(1pt);\node[below] at (v0) {${v_0}$};
        \coordinate (v1) at (1.24,-2.16);\filldraw[black] (v1) circle(1pt);\node[below] at (v1) {${v_1}$};
        \coordinate (v2) at (2.16,-1.24);\filldraw[black] (v2) circle(1pt);\node[below right] at (v2) {${v_2}$};
        \coordinate (v3) at (2.5,0);\filldraw[black] (v3) circle(1pt);\node[right] at (v3) {${v_3}$};
        \coordinate (v4) at (2.16,1.24);\filldraw[black] (v4) circle(1pt);\node[above right] at (v4) {${v_4}$};
        \coordinate (v5) at (1.24,2.16);\filldraw[black] (v5) circle(1pt);\node[above right] at (v5) {${v_5}$};
        \coordinate (v6) at (0,2.5);\filldraw[black] (v6) circle(1pt);\node[above] at (v6) {${v_6}$};
        \coordinate (v7) at (-1.24,2.16);\filldraw[black] (v7) circle(1pt);\node[above left] at (v7) {${v_7}$};
        \coordinate (v8) at (-2.16,1.24);\filldraw[black] (v8) circle(1pt);\node[above left] at (v8) {${v_8}$};
        \coordinate (v9) at (-2.5,0);\filldraw[black] (v9) circle(1pt);\node[left] at (v9) {${v_9}$};
        \coordinate (v10) at (-2.16,-1.24);\filldraw[black] (v10) circle(1pt);\node[below left] at (v10) {${v_{10}}$};
        \coordinate (v11) at (-1.24,-2.16);\filldraw[black] (v11) circle(1pt);\node[below left] at (v11) {${v_{11}}$};
          %lines

       \draw [line width=1,green](v0) -- (v1);
       \draw [line width=1,green](v0) -- (v2);\draw [line width=1,green](v1) -- (v11);
       \draw [line width=1,green](v3) -- (v11);\draw [line width=1,green](v2) -- (v10);
       \draw [line width=1,green](v3) -- (v9);\draw [line width=1,green](v4) -- (v10);
       \draw [line width=1,green](v5) -- (v9);\draw [line width=1,green](v4) -- (v8);
       \draw [line width=1,green](v5) -- (v7);\draw [line width=1,green](v6) -- (v8);
       \draw [line width=1,green](v7) -- (v6);
       \draw [line width=1,blue](v1) -- (v2);\draw [line width=1,blue](v2) -- (v3);
       \draw [line width=1,blue](v1) -- (v4);\draw [line width=1,blue](v3) -- (v0);
       \draw [line width=1,blue](v0) -- (v5);\draw [line width=1,blue](v4) -- (v11);
       \draw [line width=1,blue](v11) -- (v6);\draw [line width=1,blue](v5) -- (v10);
       \draw [line width=1,blue](v10) -- (v7);\draw [line width=1,blue](v6) -- (v9);
       \draw [line width=1,blue](v9) -- (v8);\draw [line width=1,blue](v7) -- (v8);
       \draw [line width=1,red](v3) -- (v4);
       \draw [line width=1,red](v4) -- (v2);\draw [line width=1,red](v3) -- (v5);
       \draw [line width=1,red](v5) -- (v1);\draw [line width=1,red](v2) -- (v6);
       \draw [line width=1,red](v1) -- (v7);\draw [line width=1,red](v0) -- (v6);
       \draw [line width=1,red](v0) -- (v8);\draw [line width=1,red](v11) -- (v7);
       \draw [line width=1,red](v11) -- (v9);\draw [line width=1,red](v10) -- (v8);
       \draw [line width=1,red](v10) -- (v9);
       \draw [line width=1,yellow](v4) -- (v5);
       \draw [line width=1,yellow](v5) -- (v6);\draw [line width=1,yellow](v6) -- (v3);
       \draw [line width=1,yellow](v4) -- (v7);\draw [line width=1,yellow](v7) -- (v2);
       \draw [line width=1,yellow](v3) -- (v8);\draw [line width=1,yellow](v8) -- (v1);
       \draw [line width=1,yellow](v2) -- (v9);\draw [line width=1,yellow](v9) -- (v0);
       \draw [line width=1,yellow](v1) -- (v10);\draw [line width=1,yellow](v10) -- (v11);
       \draw [line width=1,yellow](v0) -- (v11);

       \filldraw[black] (v0) circle(2pt);
      \filldraw[black] (v1) circle(2pt);\filldraw[black] (v2) circle(2pt);\filldraw[black] (v3) circle(2pt);\filldraw[black] (v4) circle(2pt);
      \filldraw[black] (v5) circle(2pt);\filldraw[black] (v6) circle(2pt);\filldraw[black] (v7) circle(2pt);\filldraw[black] (v8) circle(2pt);
      \filldraw[black] (v9) circle(2pt);\filldraw[black] (v10) circle(2pt);\filldraw[black] (v11) circle(2pt);
\end{tikzpicture}
\centering
(a)
\end{minipage}
\begin{minipage}{.45\textwidth}
\begin{tikzpicture}

              %points
       \coordinate (v0) at (0.42,-2.46);\filldraw[black] (v0) circle(2pt);\node[below right] at (v0) {${v_0}$};
       \coordinate (v1) at (1.52,-2);\filldraw[black] (v1) circle(2pt);\node[below right] at (v1) {${v_1}$};
       \coordinate (v2) at (2.28,-1.06);\filldraw[black] (v2) circle(2pt);\node[below right] at (v2) {${v_2}$};
       \coordinate (v3) at (2.5,0);\filldraw[black] (v3) circle(2pt);\node[right] at (v3) {${v_3}$};
       \coordinate (v4) at (2.2,1.18);\filldraw[black] (v4) circle(2pt);\node[above right] at (v4) {${v_4}$};
       \coordinate (v5) at (1.4,2.08);\filldraw[black] (v5) circle(2pt);\node[above right] at (v5) {${v_5}$};
       \coordinate (v6) at (0.28,2.49);\filldraw[black] (v6) circle(2pt);\node[above ] at (v6) {${v_6}$};
       \coordinate (v7) at (-0.92,2.36);\filldraw[black] (v7) circle(2pt);\node[above left] at (v7) {${v_7}$};
       \coordinate (v8) at (-1.9,1.64);\filldraw[black] (v8) circle(2pt);\node[above left] at (v8) {${v_8}$};
       \coordinate (v9) at (-2.44,0.56);\filldraw[black] (v9) circle(2pt);\node[left] at (v9) {${v_9}$};
       \coordinate (v10) at (-2.4,-0.7);\filldraw[black] (v10) circle(2pt);\node[left] at (v10) {${v_{10}}$};
       \coordinate (v11) at (-1.8,-1.74);\filldraw[black] (v11) circle(2pt);\node[below left] at (v11) {${v_{11}}$};
       \coordinate (v12) at (-0.78,-2.4);\filldraw[black] (v12) circle(2pt);\node[below] at (v12) {${v_{12}}$};
         %lines
      \draw [line width=1,green](v0) -- (v1); \draw [line width=1,green](v1) -- (v12); \draw [line width=1,green](v12) -- (v3);
      \draw [line width=1,green](v3) -- (v10); \draw [line width=1,green](v10) -- (v5); \draw [line width=1,green](v5) -- (v8);
      \draw [line width=1,green](v8) -- (v7); \draw [line width=1,green](v7) -- (v6); \draw [line width=1,green](v6) -- (v9);
      \draw [line width=1,green](v9) -- (v4); \draw [line width=1,green](v4) -- (v11); \draw [line width=1,green](v11) -- (v2);
      \draw [line width=1,green](v2) -- (v0);

      \draw [line width=1,blue](v8) -- (v9); \draw [line width=1,blue](v9) -- (v7); \draw [line width=1,blue](v7) -- (v11);
      \draw [line width=1,blue](v11) -- (v5); \draw [line width=1,blue](v5) -- (v0); \draw [line width=1,blue](v0) -- (v3);
      \draw [line width=1,blue](v3) -- (v2); \draw [line width=1,blue](v2) -- (v1); \draw [line width=1,blue](v1) -- (v4);
      \draw [line width=1,blue](v4) -- (v12); \draw [line width=1,blue](v12) -- (v6); \draw [line width=1,blue](v6) -- (v10);
      \draw [line width=1,blue](v10) -- (v8);

      \draw [line width=1,red](v3) -- (v4); \draw [line width=1,red](v4) -- (v2); \draw [line width=1,red](v2) -- (v6);
      \draw [line width=1,red](v6) -- (v0); \draw [line width=1,red](v0) -- (v8); \draw [line width=1,red](v8) -- (v11);
      \draw [line width=1,red](v11) -- (v10); \draw [line width=1,red](v10) -- (v9); \draw [line width=1,red](v9) -- (v12);
      \draw [line width=1,red](v12) -- (v7); \draw [line width=1,red](v7) -- (v1); \draw [line width=1,red](v1) -- (v5);
      \draw [line width=1,red](v5) -- (v3);

      \draw [line width=1,yellow](v11) -- (v12); \draw [line width=1,yellow](v12) -- (v10); \draw [line width=1,yellow](v10) -- (v1);
      \draw [line width=1,yellow](v1) -- (v8); \draw [line width=1,yellow](v8) -- (v3); \draw [line width=1,yellow](v3) -- (v6);
      \draw [line width=1,yellow](v6) -- (v5); \draw [line width=1,yellow](v5) -- (v4); \draw [line width=1,yellow](v4) -- (v7);
      \draw [line width=1,yellow](v7) -- (v2); \draw [line width=1,yellow](v2) -- (v9); \draw [line width=1,yellow](v9) -- (v0);
      \draw [line width=1,yellow](v0) -- (v11);

   %    \node at (v1) {\textbullet};
     \filldraw[black] (v1) circle(2pt);\filldraw[black] (v2) circle(2pt);\filldraw[black] (v3) circle(2pt);\filldraw[black] (v4) circle(2pt);
     \filldraw[black] (v5) circle(2pt);\filldraw[black] (v6) circle(2pt);\filldraw[black] (v7) circle(2pt);\filldraw[black] (v8) circle(2pt);
     \filldraw[black] (v9) circle(2pt);\filldraw[black] (v10) circle(2pt);\filldraw[black] (v11) circle(2pt);\filldraw[black] (v12) circle(2pt);
     \filldraw[black] (v0) circle(2pt);

\end{tikzpicture}
\centering
(b)
\end{minipage}
}
\caption{\small {1-PHCs on point sets in convex position: (a) $n=12\/$ and (b) $n=13\/$. }  \label{con1}}

\end{figure}

\vspace{1mm} Now, we will show how to pack $\lfloor\frac{n}{3}\rfloor$ 1-PHCs into $K_n\/$. To ensure that the boundary edges of all cycles $\{C_1, C_2, ... , C_k\}\/$ are in consecutive order, let $\{C_1, C_2, ... , C_k\}\/$ be divided into two sets $A\/$ and $B\/$ where each 1-PHC in $A\/$ has only two boundary edges, and each 1-PHC in $B\/$ has only four boundary edges (which are two couples of boundary edges and each couple has two boundary edges in consecutive order). We arrange the boundary edges with the property that a single boundary edge in $C\in A\/$ is in consecutive order with a couple of boundary edges in $C' \in B\/$ . This property is depicted in Figure \ref{con1}(a).\

\vspace{1mm} For each $i=0,1,...,\lfloor \frac{m}{3}\rfloor-1$ and $j=0,1,...,\lfloor \frac{m}{3}\rfloor-1$ where $m\geq 2\/$. Let\

$C_i = v_{3i} v_{3i+1} v_{3i-1} v_{3i+3} v_{3i-3} ... v_{3i-m} v_{2i+(m+1)},v_{3i-m-1}v_{2i+(m+3)} ... v_{3i-(2m-2)} v_{3i}\/$ and $C_j = v_{3j+2} v_{3j+1}v_{3j+4}v_{3j-1} ... v_{3j+m+1} v_{3j-m+2},v_{3j+m+3} v_{3j-m}, ... v_{3j+(2m)} v_{3j-(2m-3)}v_{3j+2}\/$.

\vspace{1mm} Here the operations on the subscripts are reduced modulo $n-1\/$.\

\vspace{1mm}(2) Now, let $n=2m+1\/$, which is odd. By Lemma \ref{lem1}, every 1-PHC in $P$ contains at least three boundary edges. On the other hand, $P\/$ has $n$ edges. Therefore, the number of 1-PHCs that can be packed into $K_n\/$ is at most $\lfloor \frac{n}{3}\rfloor\/$.\

\vspace{1mm} Now, we will show that how to pack $\lfloor \frac{n}{3}\rfloor$ 1-PHCs into $K_n\/$. To ensure that the boundary edges of all cycles $\{C_1, C_2, ... , C_k\}\/$ are in consecutive order, let each 1-PHC have two consecutive ordered boundary edges and one single boundary edge. We arrange the boundary edges with the property that a single boundary edge in $C\/$ is in consecutive order with two consecutive boundary edges in $C'\/$ and vice versa. This property is depicted in Figure \ref{con1}(b).\

\vspace{1mm} For each $i=0,1,...,\lfloor \frac{n}{3}\rfloor-1$. Let \

$C_i = v_{(m+2)i}\ v_{(m+2)i+1}\ v_{(m+2)i-1}\ v_{(m+2)i+3}\ v_{(m+2)i-3}\ ... \ v_{(m+2)i-(2m-1)}\ v_{(m+2)i}$.

\vspace{1mm} Here the operations on the subscripts are reduced modulo $n-1\/$.  \qed

\vspace{2mm}

In the next section, we will require the following additional result.

\begin{support} \label{lem2}

Let $P\/$ be a set of $n\/$ points in convex position in the plane where $n\geq 3\/$. Suppose $T\/$ is a 1-plane Hamiltonian path on $P\/$ with two pendent vertices $v_i\/$ and $v_j\/$. Then the following statements hold:\

(1) $T\/$ has at least one boundary edge when $|i-j|=1$, and each diagonal edge in $T \/$ has at least one boundary edge on each side. \

(2) $T\/$ has at least two boundary edges when $|i-j|>1$, and each diagonal edge in $T \/$ has at least one boundary edge on each side. \

\end{support}

\vspace{1mm}  \noindent
{\bf Proof:} Let $T\/$ be a 1-plane Hamiltonian path on $P\/$ with two pendent vertices $v_i\/$ and $v_j\/$. Assume that $|i-j|=1$, and let $C=T \cup \{v_iv_j\}\/$. Then $C\/$ is a 1-PHC since $v_iv_j\/$ is a boundary edge. By Lemma \ref{lem1}, $C\/$ has at least two boundary edges when $n\/$ is even and three boundary edges when $n\/$ is odd. Note that $v_iv_j\/$ is a boundary edge in $C\/$ since $|i-j|=1$. Observe that by Proposition \ref{pro2}, each diagonal edge in $C\/$ has at least one boundary edge on each side. This proves (1).\

\vspace{1mm} Assume that $|i-j|>1$. By induction on $n\/$, if $n=4\/$, the statement is trivially true. Assume that $n\geq 5\/$ and the lemma is true when $m < n\/$.\

\vspace{1mm} In the case that all the edges in $T\/$ are boundary edges, the statement holds. Hence, assume that there is a diagonal edge $v_rv_s\in E(T)\/$.\

\vspace{1mm} Let $P_1\/$ and $P_2\/$ be two sets of points of $P$ on each side of $v_rv_s\/$ that both include $v_r\/$ and $v_s\/$. Let $T_1\/$ and $T_2\/$ be the edges of $T\/$ on $P_1\/$ and $P_2\/$, respectively. It is clear that $T_i\/$ is a 1-plane Hamiltonian path on $P_i\/$ and $P_i\/$ is in convex position, $i=1,2\/$.\

\vspace{1mm} By the induction hypothesis, $T_i\/$ has at least two boundary edges. Note that $v_rv_s\/$ is a boundary edge in $T_i\/$, but it is not boundary edge in $T\/$. This proves (2).\qed

\section{Points in wheel configuration}

In this section, we turn to another special configuration. We say a set $P\/$ of $n\/$ points, is in {\em regular wheel configuration\/} if $n-1\/$ of its points are regularly spaced on a circle $C(P)\/$ with one point $x\/$ in the center of $C(P)\/$. We call $x\/$ the {\em center} of $P\/$. Note that when $n\/$ is even, $|C(P)|\/$ is odd, and since $C(P)\/$ is regularly spaced on a circle, a line passing through any two points in $C(P)\/$ does not contain $x\/$. On the other hand, when $n\/$ is odd, $|C(P)|\/$ is even, and by regularity of $C(P)\/$, $x\/$ lies on a line that passes through any two points $v_i\/$ and $v_j\/$ in $C(P)\/$ such that $|i-j|=\frac{n-1}{2}\/$. Hence, we assume $n\/$ is even. Observe that the vertices in $C(P)\/$ are the convex hull of $P\/$. An edge of the form $xv\/$ is called a {\em radial} edge, and every 1-PHC on $P\/$ contains two radial edges. \

\vspace{1mm}
\begin{support}\label{lem3}
Let $P\/$ be a set of $n\/$ points in regular wheel configuration in the plane where $n\geq 4\/$, is even. Suppose $C\/$ is a 1-PHC on $P\/$. Then $C\/$ has at least two boundary edges and any diagonal edge in $C\/$ has at least one boundary edge on each side.
\end{support}

\vspace{1mm}  \noindent
{\bf Proof:} The lemma is trivially true when n=4. hence assume $n\geq 6\/$. Suppose $x\/$ is the center of $P\/$ and $v_0v_1 \cdots v_{n-1}v_0\/$ is the cycle $C(P)\/$. Assume that $C\/$ is a 1-PHC where $xv_i\/$ and $xv_j\/$ are two radial edges of $C\/$. It is clear that $C-\{v_ix,v_jx\}\/$ is a 1-plane Hamiltonian path on $C(P)\/$. By Lemma \ref{lem2}, $C-\{v_ix,v_jx\}\/$ has at least two boundary edges except the case when $j=i+1\/$, which possibly has only one boundary edge. Furthermore, each diagonal edge in $C-\{v_ix,v_jx\}\/$ has at least one boundary edge on each side.\

\vspace{1mm} Assume that $j=i+1\/$. Suppose that $C-\{v_ix,v_jx\}\/$ has only one boundary edge $e\/$. Let $C'=C-\{v_ix,v_jx\}\cup \{v_iv_j\} \/$. Then $C'\/$ is a 1-PHC on $C(P)\/$ and has only two boundary edges $v_iv_j\/$ and $e\/$. By Proposition \ref{pro3}, $C'\/$ has the crossing edges $v_iv_{i+2}\/$ and $v_{i+1}v_{i-1}\/$; that is, $v_iv_{i+2}\/$ and $v_{i+1}v_{i-1}\/$ are edges in $C\/$. Then the radial edges $xv_i\/$ cross $v_jv_{i-1}\/$ in $C\/$ and $xv_j\/$ crosses $v_iv_{j+1}\/$ in $C\/$ (since $n\geq 6\/$), which is a contradiction since $C\/$ is a 1-PHC. Thus $C\/$ has at least two boundary edges. \qed

\vspace{1mm}
\begin{propo}\label{pro5}
Let $P\/$ be a set of $n\/$ points in regular wheel configuration in the plane where $n\geq 8\/$, is even. Suppose $C\/$ is a 1-PHC on $P\/$. If $C\/$ has only two boundary edges $v_kv_{k+1}$ for $k\in \{r,s\}\/$. Then $C\/$ has the edges $v_kv_{k+2}\/$ and $v_{k+1}v_{k-1}\/$.
\end{propo}

\vspace{1mm}  \noindent
{\bf Proof:} It is not difficult to verify that the proposition is not true when $n=6\/$. Hence, assume $n\geq 8\/$. Suppose $x\/$ is the center of $P\/$ and $v_0v_1 \cdots v_{n-2}v_0\/$ is in $C(P)\/$. Let $C\/$ be a 1-PHC on $P\/$ that contains only two boundary edges $v_kv_{k+1}\/$ for $k\in \{r,s\}\/$.\

\vspace{1mm} It is clear that $v_kv_{k+1}\/$ for $k\in \{r,s\}\/$ are single boundary edges otherwise, by Lemma \ref{lem3} any diagonal edge in $C\/$ where $v_rv_{r+1}\/$ and $v_sv_{s+1}\/$ on one side has a third boundary edge on the other sides, a contradiction. Before proceeding, we shall take note of the following observation.\

\vspace{1mm} (O1) If $v_pv_q\/$ any edge in $C\/$, then $r+1\leq p\leq s\/$ and $s+1\leq q\leq r\/$, otherwise, by Lemma \ref{lem3} $C\/$ has at least three boundary edge, a contradiction.\

\vspace{1mm} Assume on the contrary that at least one of the two edges $\{v_kv_{k+2}, v_{k+1}v_{k-1}\}\/$ is not in $C\/$ for some $k\in \{r,s\}\/$. Without loss of generality, assume that $v_rv_{r+2}\notin E(C)\/$. Suppose that $v_rv_i\/$ and $v_{r+1}v_j\in E(C)\/$, where $i\neq j\/$, we consider the following two cases.\

\vspace{1mm}
{\em Case (1):} Suppose that $x\notin \{v_i,v_j\}\/$. By (O1) $r+3\leq i\leq s\/$ and $s+1\leq j\leq r-1\/$. Hence, $v_rv_i\/$ and $v_{r+1}v_j\/$ are crossing. Thus, there is a vertex $v_t\/$ where $r+2\leq t<i\/$ such that at least one of the two edges that matches $v_t\/$ in $C\/$ crosses $v_rv_i\/$, a contradiction since $C\/$ is a 1-PHC.\

\vspace{1mm}
{\em Case (2):} Suppose that either $v_i= x\/$ or $v_j= x\/$. Without loss of generality, assume that $v_i= x\/$. Let $v_{r-1}v_l\/$ and $v_{r+1}v_j\/$ be two edges in $C\/$. By (O1) $s+1\leq j\leq r-1\/$ and $r+2\leq l\leq s\/$ where $v_sv_{s+1}\/$ is the second boundary edge. By regularity of $C(P)\/$, $v_rx\/$ and $v_{r+1}v_p\/$ are crossing, then $j=r-1\/$; otherwise, there is a vertex $v_t\/$ where $j+1\leq t<r$ such that at least one of the two edges that matches $v_t\/$ into $C\/$ crosses $v_{r+1}v_j\/$, in both cases there is a contradiction since $C\/$ is a 1-PHC has only two boundary edges. Thus $v_{r+1}v_{r-1}\/$ in $C\/$.\

\vspace{1mm} Now, if $l>r+3\/$. Then the edges in $C\/$ that incident on $v_{r+2},v_{r+3}\/$ crosses $v_{r-1}v_l\/$, a contradiction since $C\/$ is a 1-PHC. By regularity of $C(P)$ and $n\geq 8\/$ if $l=r+3\/$, then $v_{r-1}v_l\/$ crosses $v_rx\/$ which is a contradiction since $v_{r+1}v_{r-1} \in E(C)\/$ and crosses $v_rx\/$. This completes the proof. \qed

\vspace{3mm}
We now present the main result of this section.

\vspace{1mm}
\begin{result}\label{the2}
Let $P\/$ be a set of $n\/$ points in regular wheel configuration in the plane where $n\geq 10\/$, is even. Then there exist $k\/$ edges-disjoint 1-PHCs $C_1, C_2, ... , C_k\/$ on $P\/$ that can be packed into $K_n\/$ where $k\leq \lfloor \frac{n-1}{3}\rfloor\/$.
\end{result}

\vspace{1mm}  \noindent
{\bf Proof:} Suppose $x\/$ is the center of $P\/$ and $v_0v_1 \cdots v_{n-1}v_0\/$ is the cycle $C(P)\/$. By Lemma \ref{lem3}, every 1-PHC in $P\/$ contains at least two boundary edges. On the other hand, $C(P)\/$ has $n-1\/$ boundary edges; that is, the number of 1-PHCs does not exceed $n-1/2\/$.\

\vspace{1mm} By Proposition \ref{pro5}, if $C_i\/$ and $C_j\/$ are two edge-disjoint 1-PHCs each having only two boundary edges, then any boundary edge of $C_i\/$ can not be in consecutive order with a boundary edge of $C_j\/$.\

\vspace{1mm} Note that a single boundary edge in $C_i\/$ can be adjacent to two consecutive boundary edges in $C_j\/$ where $i\neq j\/$; that is, $C_i \cup C_j\/$ can have three boundary edges in consecutive order. Therefore, the number of 1-PHCs that can be packed into $K_n\/$ is at most $\lfloor \frac{n-1}{3}\rfloor\/$ and this bound is tight.\

\begin{figure}[htb]
\centering
\resizebox{6cm}{!}{
\begin{tikzpicture}
       %points
       \coordinate (x) at (0,0);\node at (x) {\textbullet};\node[above left] at (x) {${x}$};
       \coordinate (v0) at (0.3,-2.5);\node[below] at (v0) {${v_0}$};
       \coordinate (v6) at (0.3,2.5);\node[above] at (v6) {${v_6}$};
       \coordinate (v1) at (1.44,-2.06);\node[below] at (v1) {${v_1}$};
       \coordinate (v5) at (1.44,2.06);\node[above right] at (v5) {${v_5}$};
       \coordinate (v2) at (2.22,-1.18);\node[ right] at (v2) {${v_2}$};
       \coordinate (v4) at (2.22,1.18);\node[above right] at (v4) {${v_4}$};
       \coordinate (v3) at (2.5,0);\node[right] at (v3) {${v_3}$};
       \coordinate (v7) at (-0.9,2.34);\node[above] at (v7) {${v_7}$};
       \coordinate (v12) at (-0.9,-2.34);\node[below left] at (v12) {${v_{12}}$};
       \coordinate (v9) at (-2.44,0.6);\node[left] at (v9) {${v_{9}}$};
       \coordinate (v10) at (-2.44,-0.6);\node[left] at (v10) {${v_{10}}$};
       \coordinate (v8) at (-1.88,1.65);\node[above left] at (v8) {${v_{8}}$};
       \coordinate (v11) at (-1.88,-1.65);\node[below left] at (v11) {${v_{11}}$};
         %lines
       \draw [line width=1,green](v0) -- (v1);\draw [line width=1,green](v1) -- (v12);\draw [line width=1,green](v12) -- (v3);
       \draw [line width=1,green](v3) -- (v10);\draw [line width=1,green](v10) -- (v5);\draw [line width=1,green](v5) -- (v8);
       \draw [line width=1,green](v8) -- (v7);\draw [line width=1,green](v7) -- (v6);\draw [line width=1,green](v6) -- (v9);
       \draw [line width=1,green](v9) -- (v4); \draw [line width=1,green](v4) -- (x); \draw [line width=1,green](x) -- (v11);
       \draw [line width=1,green](v11) -- (v2);\draw [line width=1,green](v2) -- (v0);

       \draw [line width=1,blue](v8) -- (v9);\draw [line width=1,blue](v9) -- (v7);\draw [line width=1,blue](v7) -- (v11);
       \draw [line width=1,blue](v11) -- (v5);\draw [line width=1,blue](v5) -- (v0);\draw [line width=1,blue](v0) -- (v3);
       \draw [line width=1,blue](v3) -- (v2);\draw [line width=1,blue](v2) -- (v1);\draw [line width=1,blue](v1) -- (v4);
       \draw [line width=1,blue](v4) -- (v12); \draw [line width=1,blue](v12) -- (x); \draw [line width=1,blue](x) -- (v6);
       \draw [line width=1,blue](v6) -- (v10);\draw [line width=1,blue](v10) -- (v8);

       \draw [line width=1,red](v3) -- (v4);\draw [line width=1,red](v4) -- (v2);\draw [line width=1,red](v2) -- (v6);
       \draw [line width=1,red](v6) -- (v0);\draw [line width=1,red](v0) -- (v8);\draw [line width=1,red](v8) -- (v11);
       \draw [line width=1,red](v11) -- (v10);\draw [line width=1,red](v10) -- (v9);\draw [line width=1,red](v9) -- (v12);
       \draw [line width=1,red](v12) -- (v7); \draw [line width=1,red](v7) -- (x); \draw [line width=1,red](x) -- (v1);
       \draw [line width=1,red](v1) -- (v5);\draw [line width=1,red](v5) -- (v3);

       \draw [line width=1,yellow](v11) -- (v12);\draw [line width=1,yellow](v12) -- (v10);\draw [line width=1,yellow](v10) -- (v1);
       \draw [line width=1,yellow](v1) -- (v8);\draw [line width=1,yellow](v8) -- (v3);\draw [line width=1,yellow](v3) -- (v6);
       \draw [line width=1,yellow](v6) -- (v5);\draw [line width=1,yellow](v5) -- (v4);\draw [line width=1,yellow](v4) -- (v7);
       \draw [line width=1,yellow](v7) -- (v2); \draw [line width=1,yellow](v2) -- (x); \draw [line width=1,yellow](x) -- (v9);
       \draw [line width=1,yellow](v9) -- (v0);\draw [line width=1,yellow](v0) -- (v11);

       \node at (x) {\textbullet}; \filldraw[black] (v0) circle(2pt);\filldraw[black] (v1) circle(2pt);\filldraw[black] (v2) circle(2pt);
       \filldraw[black] (v3) circle(2pt);\filldraw[black] (v4) circle(2pt);\filldraw[black] (v5) circle(2pt);\filldraw[black] (v6) circle(2pt);
       \filldraw[black] (v7) circle(2pt);\filldraw[black] (v8) circle(2pt);\filldraw[black] (v9) circle(2pt);
       \filldraw[black] (v10) circle(2pt);\filldraw[black] (v11) circle(2pt);\filldraw[black] (v12) circle(2pt);

\end{tikzpicture}
}
\caption{\small {1-PHCs on a set of points in regular wheel configuration, $n=14\/$. } \label{reg1}}
\end{figure}
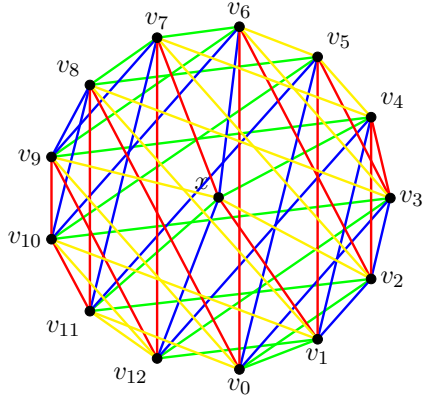

\vspace{1mm} Now, we will show how to pack $\lfloor\frac{n-1}{3}\rfloor$ 1-PHCs into  $K_n\/$. To ensure that the boundary edges of all cycles $\{C_1, C_2, ... , C_k\}\/$ are in consecutive order. Then each 1-PHC should have three boundary edges (where two of them are in consecutive order) with the property that a single boundary edge in $C_i\/$ is adjacent to the two consecutive boundary edges in $C_j\/$ and vice versa. This property is depicted in Figure  \ref{reg1}. For each $i=0,1,...,\lfloor \frac{n-1}{3}\rfloor-1$ and $r=\lfloor \frac{m+1}{2}\rfloor\/$. Let

$C_i = v_{(m+1)i} v_{(m+1)i+1} v_{(m+1)i-1} v_{(m+1)i+3} v_{(m+1)i-3} \dots  v_{(m+1)i+r} \ x \ v_{(m+1)i+5} v_{(m+1)i-5}\/$

$ \dots v_{(m+1) i-(2m-3)} v_{(m+1)i}\/$. Here the operations on the subscripts are reduced modulo $2n-1\/$. \qed

\vspace{1mm}
\section{1-PHCs for point sets in general position}
In this section, we consider a set $P\/$ of $n\/$ points in general position in the plane i.e., no three points are collinear. For $n =2^k +h\/$ where $0 \leq h < 2^k\/$, we will show that there are at least $k-1\/$ edge-disjoint 1-PHCs on $P\/$ (Theorem \ref{the3}). For this purpose, we present some ingredients that will be used to prove the main result in this section.

\subsection{Bisect lines for a set of points}
Let $P\/$ be a set of $n\/$ points in general position in the plane. A line $l\/$ is said to \textit{bisect} a set $P\/$ if both open half spaces defined by $l\/$ contain precisely $\frac{n}{2}\/$ points. It is no loss of generality to assume $n\/$ is odd since otherwise, any point $v\/$ may be removed and any line that bisects $P-\{w\}$ also bisects $P\/$.\

Let $P_1\/$ and $P_2\/$ be two point sets in the plane. If $H_1\/$ and $H_2\/$ are two convex polygons containing $P_1\/$ and $P_2\/$ respectively, then we say that $H_1\/$ and $H_2\/$ are disjoint if there is a line that separates them. Moreover, if $P\/$ is a disjoint union of two point sets $P_1\/$ and $P_2\/$, the ham-sandwich cut theorem guarantees the existence of a line that simultaneously bisects $P_1\/$ and $P_2\/$ (see for example \cite {h:refer},\cite {cmw:refer}).\

\vspace{1mm}
\begin{support} \label{lem4}
Let $P\/$ be a set of $m\/$ points in the general position where $m \geq 3\/$. Suppose there is a line separating a given set  $\{v, w\}\/$  from $P\/$. Then there is a line that separates some subset $P_1\/$ from $P\/$ with $\{v, w\} \subseteq P_1\/$ and $2 \leq |P_1| \leq m-1\/$.
\end{support}

\vspace{1mm}  \noindent
{\bf Proof:} The lemma is trivially true if $m=3\/$ with $P_1 = \{v, w\}\/$. So assume that $m \geq 4\/$.

\vspace{1mm}
Let $v_1 \in P - \{v, w\}\/$ be such that all points in $P - \{v, w, v_1\}\/$ are on one side of the line $l_1\/$ joining $v_1\/$ and $z_1\/$ for some $z_1 \in \{v, w\}\/$ and let $P_1 = \{v, w, v_1\}\/$. Let $L_1\/$ be a line parallel to $l_1\/$ such that all points in $P - \{v, w, v_1\}\/$ are on one side of $L_1\/$.

\vspace{1mm} If $|P_1| = m-1\/$, then the proof is complete.  Otherwise  repeat the argument with $v_2 \in P - \{v, w, v_1\}\/$ and $z_2 \in  \{v, w, v_1\}\/$ so that all points in $P - \{v, w, v_1, v_2\}\/$ are on one side of the line $l_2\/$ joining $v_2z_2\/$ and let $P_1 = \{v, w, v_1, v_2\}\/$ with the line $L_2\/$ similarly defined. By repeating the argument where necessary, we reach the conclusion of the lemma. \qed

\vspace{1mm}
\begin{support} \label{lem5}
Let $L\/$ be a line bisects a set $P\/$ of $m\/$ points in the general position into $P_1\/$ and $P_2\/$ where $m \geq 6\/$, and let $l^{\bot}\/$ be a line perpendicular to $L\/$ and all points in $P\/$ are on one side of $l^{\bot}\/$. Suppose $\{u, v\}\subset P_1\/$ is a given set such that no more than half points of $P_1\/$ are between $l^{\bot}\/$ and a line, perpendicular to $L\/$, passing through any point in $\{v, w\}\/$. Then there is a line bisects $P_1\/$ and $P_2\/$ into $P_{i,j}\/$, for each $i=1,2\/$ with $j=1,2\/$ and $\{v, w\} \subseteq P_{1,k}\/$ for some $k\in\{1,2\}\/$.

\end{support}

\vspace{1mm}  \noindent
{\bf Proof:} The lemma is trivially not true if $m=5\/$ with $P_1 = \{v, w\}\/$. So assume that $m \geq 6\/$. By hum-sandwich cut theorem, there is a line $l\/$ bisects $P_1\/$ and $P_2\/$ in the plane into sets $P_{i,j}\/$, for each $i=1,2\/$ with $j=1,2\/$. It is no loss of generality to assume $P_{1,2}\/$ is between $l\/$ and $l^{\bot}\/$. Assume on the contrary that $\{v, w\} \nsubseteq P_{1,j}\/$ for each $j=1,2\/$. Then all points of $P_{1,2}\/$ are between $l^{\bot}\/$ and the line, perpendicular to $L\/$, passing through the point in $\{v, w\}\cap P_{1,1}\/$, a contradiction. Thus $\{v, w\} \subseteq P_{1,2}\/$. \qed

\subsection{Drawing a 1-PHC on a set of points}

We shall give a description of an algorithm for drawing a 1-PHC on a set of points in general position in the plane. In what follows, let $l(v_1, v_2)$ be a line passing through the two points $v_1$ and $v_2$. \

\vspace{8mm}
\textbf{Algorithm $A\/$}
\begin{enumerate}

\item Find a line $l\/$ that bisects $P\/$ into $P_1\/$ and $P_2\/$ where either $|P_1| = |P_2|\/$ or $|P_1| = |P_2| +1\/$.

\item Find a line $l^{\bot}\/$ such that $l^{\bot}\/$ is perpendicular to $l\/$ and all points in $P\/$ are on one side of  $l^{\bot}\/$.

\item Find $CH(P_i)\/$, the convex hull of $P_i\/$, for each $i=1, 2\/$ and select $v_i \in CH(P_i)\/$ such that all the points in $P_1 \cup P_2 -\{v_1, v_2\}\/$ are between $l^{\bot}\/$ and the line $l(v_1, v_2)$. Let $v_1v_2\/$ be an edge in $C\/$ and let $v^{*}_i=v_i$, for each $i=1,2\/$.

\item If $P_i = \{v^*_i\}\/$ for $i=1, 2$, let $v^{*}_1v^{*}_2\/$ be an edge in $C\/$ and Stop. If $P_2=\{v^*_2\}\/$ and $P_1 = \{v^*_1, w\}\/$ let $v^{*}_1w\/$ and $v^{*}_2w\/$ be edges in $C\/$ and Stop. Otherwise, let $P_i=P_i-\{v^*_i\}$ for each $i=1, 2$.

\item Find $CH(P_i)\/$, for $i=1, 2\/$ and select $v_i \in CH(P_i)\/$, $i=1, 2\/$, be such that all the points in $P_1 \cup P_2 -\{v^*_1, v^*_2\}\/$ are between $l^{\bot}\/$ and the line $l(v_1, v_2)$.

\item If no point of $\{v^{*}_1,v^{*}_2\}\/$, is between $l^{\bot}\/$ and the line $l(v_1, v_2)$. Let $v^{*}_1v_2\/$ and $v^{*}_2v_1\/$ be edges in $C\/$.  Repeat Step (4) with $v_i\/$ taking the place of $v^{*}_i$ for each $i=1,2\/$.

\item For some $i\in\{1,2\}\/$, if $v^{*}_{3-i}\/$ is not between $l^{\bot}\/$ and the line $l(v_i, v^{*}_i)$ and all points in $P_1 \cup P_2 -\{v^{*}_{3-i}\}\/$ are between $l^{\bot}\/$ and the line $l(v_i, v^{*}_i)$, let $v_iv^{*}_{3-i}\/$ be two edges in $C\/$.

\item If $\{v_iv^{*}_{3-i}, v_iv_{3-i}\}\subset E(C)\/$, let $P_i=P_i-\{v_i\}$ and repeat Step (4) with $v_{3-i}\/$ taking the place of $v^{*}_{3-i}$.

\end{enumerate}

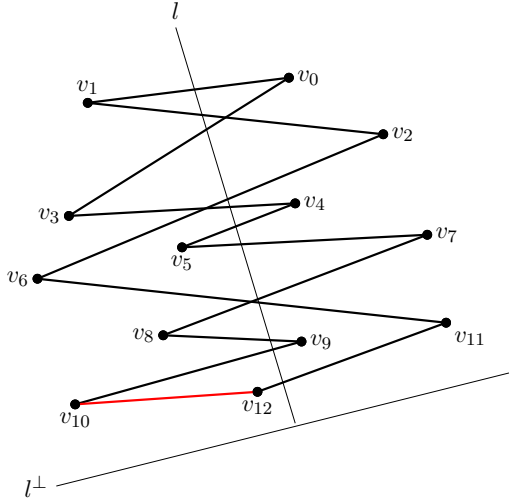
\begin{figure}[htp]
\centering
\resizebox{7cm}{!}{
\begin{tikzpicture}%[rotate=-11,.style={draw}]
%\coordinate (center) at (0,0);
 % \def\radius{2.5cm}
  % \draw (center) circle[radius=\radius];
   %\foreach \x in {0,14.4,...,360} {
    %         \filldraw[] (\x:2.5cm) circle(1pt);
     %        \filldraw[] (0,0) circle(1pt);
      %         }
      %   \draw [line width=1,green](0.22,-3.5) -- (0,0);\draw [line width=1,blue](0,0) -- (1.08,-3.3);\draw [line width=1,red](0,0) -- (1.9,-3);
           %points

         \coordinate (v0) at (1.5,1.5); \filldraw[black] (v0) circle(2pt);\node[right] at (v0) {${v_0}$};
         \coordinate (v1) at (-1.7,1.1); \filldraw[black] (v1) circle(2pt);\node[above] at (v1) {${v_{1}}$};
         \coordinate (v2) at (3,0.6); \filldraw[black] (v2) circle(2pt);\node[ right] at (v2) {${v_2}$};
         \coordinate (v3) at (-2,-0.7); \filldraw[black] (v3) circle(2pt);\node[left] at (v3) {${v_3}$};
         \coordinate (v4) at (1.6,-0.5);\filldraw[black] (v4) circle(2pt);\node[right] at (v4) {${v_{4}}$};
         \coordinate (v5) at (-0.2,-1.2); \filldraw[black] (v5) circle(2pt);\node[below] at (v5) {${v_5}$};
         \coordinate (v6) at (-2.5,-1.7);\filldraw[black] (v6) circle(2pt);\node[left] at (v6) {${v_{6}}$};
         \coordinate (v7) at (3.7,-1); \filldraw[black] (v7) circle(2pt);\node[right] at (v7) {${v_{7}}$};
         \coordinate (v8) at (-0.5,-2.6); \filldraw[black] (v8) circle(2pt);\node[left] at (v8) {${v_{8}}$};
         \coordinate (v9) at (1.7,-2.7);\filldraw[black] (v9) circle(2pt);\node[right] at (v9) {${v_{9}}$};
         \coordinate (v10) at (-1.9,-3.7);\filldraw[black] (v10) circle(2pt);\node[below] at (v10) {${v_{10}}$};
         \coordinate (v11) at (4,-2.4);\filldraw[black] (v11) circle(2pt);\node[below right] at (v11) {${v_{11}}$};
         \coordinate (v12) at (1,-3.5);\filldraw[black] (v12) circle(2pt);\node[below] at (v12) {${v_{12}}$};
       %%%%%%%%%%%%%%%%%%%%%%%%%%%%%%%%%%%%%%%%%%%%%%%%%%%%%%%%%%%%%%%%%%%%%%%%%%%%%%%%%%%%%%%%%%%%%%%%%%%%%%%%%%%%%%%%
       \coordinate (v25) at (1.6,-4);\coordinate (v26) at (-0.3,2.3);\node[above] at (v26) {${l}$};
       \draw [line width=.1,black](v25) -- (v26);
       \coordinate (v29) at (5,-3.2);\coordinate (v30) at (-2.2,-5);\node[left] at (v30) {$l^{\bot}\/$};
       \draw [line width=.1,black](v29) -- (v30);

        %%%%%%%%%%%%%%%%%%%%%%%%%%%%%%%%%%%%%%%%%%%%%%%%%%%%%%%%%%%%%%%%%%%%%%%%%%%%%%%%%%%%%%%%%%%%%%%%%%%%%%%%%%%%%%%%
         %lines
     \draw [line width=1,black](v1) -- (v0);
     \draw [line width=1,black](v1) -- (v2);\draw [line width=1,black](v0) -- (v3);
     \draw [line width=1,black](v2) -- (v6);\draw [line width=1,black](v3) -- (v4);
     \draw [line width=1,black](v4) -- (v5);\draw [line width=1,black](v5) -- (v7);
     \draw [line width=1,black](v6) -- (v11);\draw [line width=1,black](v7) -- (v8);
     \draw [line width=1,black](v8) -- (v9);\draw [line width=1,black](v11) -- (v12);
     \draw [line width=1,black](v9) -- (v10);\draw [line width=1,red](v10) -- (v12);
     %%%%%%%%%%%%%%%%%%%%%%%%%%%%%%%%%%%%%%%%%%%%%%%%%%%%%%%%%%%%%%%%%%%%%%%%%%%%%%%%%%%%%%%%%%%%%%%%%%%%%%%%%%%%%%%

  %    \node at (v1) {\textbullet};
    \filldraw[black] (v0) circle(2pt);\filldraw[black] (v1) circle(2pt);\filldraw[black] (v2) circle(2pt);\filldraw[black] (v3) circle(2pt);
     \filldraw[black] (v4) circle(2pt);\filldraw[black] (v5) circle(2pt);\filldraw[black] (v6) circle(2pt);\filldraw[black] (v7) circle(2pt);
     \filldraw[black] (v8) circle(2pt);\filldraw[black] (v9) circle(2pt);\filldraw[black] (v10) circle(2pt);\filldraw[black] (v11) circle(2pt);
     \filldraw[black] (v12) circle(2pt);
\end{tikzpicture}
}
\caption{\small Drawing a 1-PHC on a set of $13\/$ points in general position by Algorithm $A\/$, and the red edge $v_{10}v_{12}\/$ refer to the Stone $st(v_{10},v_{12})\/$.}  \label{algorithm}
\end{figure}

\vspace{1mm}
The edge $v^{*}_1w\/$ in Step (4) is termed a "stone" and shall be denoted by $st(v, w)\/$. A 1-PHC obtained by Algorithm (A) which contains a stone, is depicted in Figure \ref{algorithm}.

\vspace{1mm}
\subsection{A joining between two 1-PHCs }
In this section, we show how to extract a 1-PHC by joining two edge-disjoint 1-PHCs. Let $P(1)\/$ and $P(2)\/$ be two disjoint point sets in general position in the plane. Suppose $C(1)\/$ and $C(2)\/$ are two edge-disjoint 1-PHCs on $P(1)\/$ and $P(2)\/$, respectively.\

The edges $u_1u_2\in E(C(1))\/$ and $v_1v_2\in E(C(2))\/$ are called \emph{joining edges} of $C(1)\/$ and $C(2)\/$ if the graph resulting from removing them and adding the edges $u_1v_1\/$ and $u_2v_2\/$ (or $u_1v_2\/$ and $u_2v_1\/$) still an edge-disjoint 1-PHC on $P_1\cup P_2\/$. The edges $u_1v_1\/$ and $u_2v_2\/$ (or $u_1v_2\/$ and $u_2v_1\/$) are termed a \emph{connection edges} of $C(1)\/$ and $C(2)\/$.\

Suppose that there are two crossing edges $u_1u_2\/$ and $u_3u_4\/$ in $C(1)\/$ such that the graph resulting from removing them and adding two non-crossing edges $u_1u_3\/$ and $u_2u_4\/$ still a 1-PHC. Then the edges $u_1u_3\/$ and $u_2u_4\/$ are termed \emph{created} edges in $C(1)\/$. Two joining edges of $C(1)\/$ and $C(2)\/$ are called  \emph{created joining} edges if at least one of them is a created edge.\

\begin{support} \label{lem6}
Let $P\/$ be a set of $n\/$ points in general position in the plane where $n\geq 8\/$, and let $C\/$ be a 1-PHC on $P\/$ which is bisected into two disjoint sets $P_1\/$ and $P_2\/$. Suppose $C(1)\/$ and $C(2)\/$ are two 1-PHCs on $P_1\/$ and $P_2\/$, respectively, where $C(1)\/$ and $C(2)\/$ have no joining edges. Then $C(1)\/$ and $C(2)\/$ can be joined by two joining created edges.
\end{support}

\vspace{1mm}  \noindent
{\bf Proof:} Let $C\/$ be a 1-PHC on $P\/$. In so doing, the set $P\/$ has been split into $P_1\/$ and $P_2\/$. Assume that $C(1)\/$ and $C(2)\/$ are two 1-PHCs on $P_1\/$ and $P_2\/$, respectively since $n\geq 8\/$, and let $C(1)\/$ and $C(2)\/$ have no joining edges.\

\vspace{1mm}
{\em Case (1):} When $C(i)\/$ has at least one crossing for each $i=1, 2$ (since $n\geq 8\/$ and $|C(1)|\geq 4\/$ and $|C(2)|\geq 4\/$). Let $\{u_1u_2,u_3u_4\}\/$ and $\{v_1v_2,v_3v_4\}\/$ be the two crossing edges in $C(1)\/$ and $C(2)\/$, respectively. By removing the crossing edges in $C(1)\/$ and $C(2)\/$ and adding the non-crossing edges $\{u_1u_4,u_2u_3\}\/$ and $\{v_1v_4,v_2v_3\}\/$ in $C(1)\/$ and $C(2)\/$, respectively we obtain two created edges in each of $C(1)\/$ and $C(2)\/$. Chose $e_1\in\{u_1u_4,u_2u_3\}$ and $e_2\in\{v_1v_4,v_2v_3\}\/$ such that no points between them. Without loss of generality assume that no points between the two created edges $u_1u_4\/$ and $v_1v_4\/$.\

\vspace{1mm} Note that at least one edge in a set $A=\{u_1v_1,u_1v_4,u_4v_1,u_4v_4\}\/$ is not in $C\/$ since $A\/$ is 4-cycle graph and $|V(C)|\geq8\/$ is not union of two cycles. That means\

(1) If $C\/$ has only three edges of $A\/$. Then there is an edge $u'v' \in A\cap E(C)\/$ such that $d_C(u')=d_C(v')=2\/$. Thus there are two joining edges one in original $C(1)\/$ incident on $u'\/$ and another in original $C(2)\/$ incident on $v'\/$, a contradiction (since $C(1)\/$ and $C(2)\/$ have no joining edges).\

(2) If $C\/$ has only two edges $\{uv,u'v'\}\subset A\/$. Then (i) $\{uv,u'v'\}$ share no vertex, and hence $u_1u_4\/$ and $v_1v_4\/$ are joining created edges since $A-\{uv,u'v'\}\/$ are connection edges. (ii) $\{uv,u'v'\}$ share on a vertex. It is no loss of generality to assume that $u=u'=u_1\/$. That is,  $\{u_1v_1,u_1v_4\}\subset A\/$, then the created edge $u_1v_1\/$ and the edge incident on $u_1\/$ such as $u_1u^*\/$ for some $u^*\notin A\/$ are created joining edges since $A-\{u_1v_1,u_1v_4\}\/$ is not in $C\/$ and at most one of the two edges of $u^*v_1,u^*v_4 \/$ is not in $C\/$ since $A-\{u_1v_1,u_1v_4\}\cup \{u^*v_1,u^*v_4\}\/$ is 4-cyles. \

(O3) If $C\/$ has only one edge $u'v'\in A\/$. Then $u_1u_4\/$ and $v_1v_4\/$ are joining created edges since there are two connections edge in $A-\{u'v'\}\/$.

\vspace{1mm}
{\em Case (2):} When $C(i)\/$ has at most one crossing for some $i\in\{1, 2\}$. By removing the crossing edges and adding two created edges we obtain two plane cycles $C(1)\/$ and $C(2)\/$. As in case(1), asume that no points between the two edges $u_1u_2 \in E(C(1))\/$ and another $v_1v_2 \in E(C(2))\/$. By repeat the similar argument in case (1) we see $C(1)\/$ and $C(2)\/$ have two joining created edges. \qed

\subsection{Packing 1-PHCs into a point set }
We conclude this paper with the following main result.\

\vspace{1mm}
\begin{result} \label{the3}
Let $P\/$ be a set of $n\/$ points in general position in the plane where $n=2^k+h\/$, with $0\leq h <2^k\/$. Then there exist at least $k-1\/$ edge-disjoint 1-PHCs $C_1, C_2, ... , C_{k-1}\/$ on $P\/$ that can be packed into $K_n\/$.
\end{result}

\vspace{1mm}  \noindent
{\bf Proof:}  First we apply Algorithm $(A)\/$ to obtain the first 1-PHC $C_1\/$. In so doing, the set $P\/$ has been bisected into $P_1\/$ and $P_2\/$ by $l_1$. Let $P_1$ on the left of $l_1$ and $P_2$ on the right of $l_1$.

\vspace{1mm} If $P_1\/$ has no stone, then by hum-sandwich cut theorem there is a line $l_2\/$ that simultaneously bisects $P_1\/$ and $P_2\/$ into $P_{i,j}\/$, $i=1,2\/$ with $j=1,2\/$  which are label in the anticlockwise order and either $|P_{1, 1}| = |P_{1, 2}|\/$ or  $|P_{1, 1}| = |P_{1, 2}|+1\/$. \

\vspace{1mm} If $P_1\/$ has a stone $st(v, w)\/$, then we have two cases \

\vspace{1mm}
{\em Case (1):} By Lemma \ref{lem5}, there is a line $l_2\/$ that simultaneously bisects $P_1\/$ and $P_2\/$ into $P_{i,j}\/$, $i=1,2\/$ with $j=1,2\/$. and  $\{v, w\} \subsetneq P_{1,2}\/$ and either $|P_{1, 1}| = |P_{1, 2}|\/$ or  $|P_{1, 1}| = |P_{1, 2}|+1\/$.\

\vspace{1mm}
{\em Case (2):} By Lemma \ref{lem4}, there is a line $l_2\/$ that bisects $P_1\/$ into $P_{1,1}\/$ and $P_{1,2}\/$ and with $\{v, w\} \subsetneq P_{1,2}\/$ and either $|P_{1, 1}| = |P_{1, 2}|\/$ or  $|P_{1, 1}| = |P_{1, 2}|+1\/$. Furthermore, there is a line $l_2'\/$ that bisects $P_2\/$ into $P_{2,1}\/$ and $P_{2,2}\/$ and either $|P_{2, 1}| = |P_{2, 2}|\/$ or  $|P_{2, 1}| = |P_{2, 2}|+1\/$.

\vspace{1mm}In all cases, label the parts $P_{i,j}\/$ in the anticlockwise order. In case (1) and case (2), $C(1)\/$ and $C(2)\/$ are two edge-disjoint cycles can be joined using either the joining edges or created joining edges (by Lemma \ref{lem6}).\

\vspace{1mm} To obtain $C_3\/$, rename $P_{i,j}$ to be four parts $P_1$, $P_2$, $P_3$, and $P_4$ arranged in anticlockwise ordered.
 Then repeat the above operations with $P_i$ taking place $P$ for each $i=1, 2, 3, 4$ to obtain 1-PHCs $C(i)$.
  Join $C(i)$ with $C(i+1)$ for $i=1, 2, 3$ either by joining edges or by created joining edges.

\vspace{3mm} In general, to obtain $C_r\/$ where $1\leq r\leq k-1\/$,  repeat the above operations on parts $P_1, P_2, \ldots ,P_{2^{r-1}}$, arranged in anticlockwise ordered, with $P_i$ taking place $P$ for each $i=1, 2, \ldots , 2^{r-1}$ to obtain 1-PHCs $C(i)$.
  Join $C(i)$ with $C(i+1)$ for $i=1, 2, \ldots , 2^{r-1}-1$ either by joining edges or by created joining edges. \qed

\begin{figure}[htp]
\centering
\resizebox{12cm}{!}{
\begin{tikzpicture}%[rotate=-11,.style={draw}]
%\coordinate (center) at (0,0);
 % \def\radius{2.5cm}
  % \draw (center) circle[radius=\radius];
   %\foreach \x in {0,14.4,...,360} {
    %         \filldraw[] (\x:2.5cm) circle(1pt);
     %        \filldraw[] (0,0) circle(1pt);
      %         }
      %   \draw [line width=1,green](0.22,-3.5) -- (0,0);\draw [line width=1,blue](0,0) -- (1.08,-3.3);\draw [line width=1,red](0,0) -- (1.9,-3);
           %points
         \coordinate (v0) at (1.9,3); \filldraw[black] (v0) circle(2pt);\node[above right] at (v0) {${v_{0}}$};
         \coordinate (v1) at (-0.54,2.7); \filldraw[black] (v1) circle(2pt);\node[above] at (v1) {${v_{1}}$};
         \coordinate (v2) at (-3,1.5);\filldraw[black] (v2) circle(2pt);\node[above] at (v2) {${v_{2}}$};
         \coordinate (v3) at (2,1.5);\filldraw[black] (v3) circle(2pt);\node[above] at (v3) {${v_{3}}$};
         \coordinate (v4) at (0.5,0.8); \filldraw[black] (v4) circle(2pt);\node[right] at (v4) {${v_4}$};
         \coordinate (v5) at (-1.7,0.5); \filldraw[black] (v5) circle(2pt);\node[above] at (v5) {${v_{5}}$};
         \coordinate (v6) at (3.8,0.6); \filldraw[black] (v6) circle(2pt);\node[ right] at (v6) {${v_6}$};
         \coordinate (v7) at (-3,-0.7); \filldraw[black] (v7) circle(2pt);\node[left] at (v7) {${v_7}$};
         \coordinate (v8) at (1.6,-0.5);\filldraw[black] (v8) circle(2pt);\node[below left] at (v8) {${v_{8}}$};
         \coordinate (v9) at (-0.2,-1.2); \filldraw[black] (v9) circle(2pt);\node[below] at (v9) {${v_9}$};
         \coordinate (v10) at (-2.5,-1.7);\filldraw[black] (v10) circle(2pt);\node[left] at (v10) {${v_{10}}$};
         \coordinate (v11) at (3.7,-1); \filldraw[black] (v11) circle(2pt);\node[right] at (v11) {${v_{11}}$};
         \coordinate (v12) at (-0.5,-2.6); \filldraw[black] (v12) circle(2pt);\node[left] at (v12) {${v_{12}}$};
         \coordinate (v13) at (1.7,-2.7);\filldraw[black] (v13) circle(2pt);\node[right] at (v13) {${v_{13}}$};
         \coordinate (v14) at (-1.9,-3.7);\filldraw[black] (v14) circle(2pt);\node[below] at (v14) {${v_{14}}$};
         \coordinate (v15) at (4,-2.4);\filldraw[black] (v15) circle(2pt);\node[below right] at (v15) {${v_{15}}$};
         \coordinate (v16) at (1,-3.5);\filldraw[black] (v16) circle(2pt);\node[below] at (v16) {${v_{16}}$};
       %%%%%%%%%%%%%%%%%%%%%%%%%%%%%%%%%%%%%%%%%%%%%%%%%%%%%%%%%%%%%%%%%%%%%%%%%%%%%%%%%%%%%%%%%%%%%%%%%%%%%%%%%%%%%%%%
       \coordinate (v25) at (1.6,-5);
       \coordinate (v26) at (-0.3,3.6);\node[above] at (v26) {${l_1}$};
       \draw [line width=.1,black](v25) -- (v26);
       \coordinate (v27) at (5.5,1.05);
       \coordinate (v28) at (-5,-2.8);\node[left] at (v28) {${l_2}$};
        \draw [line width=.1,black](v27) -- (v28);
         \coordinate (v29) at (6,-1.2);
       \coordinate (v30) at (-4.5,-4);\node[left] at (v30) {${l_3}$};
        \draw [line width=.1,black](v29) -- (v30);
         \coordinate (v31) at (5,2.2);
       \coordinate (v32) at (-5,-0.68);\node[left] at (v32) {${l_4}$};
        \draw [line width=.1,black](v32) -- (v31);
        %%%%%%%%%%%%%%%%%%%%%%%%%%%%%%%%%%%%%%%%%%%%%%%%%%%%%%%%%%%%%%%%%%%%%%%%%%%%%%%%%%%%%%%%%%%%%%%%%%%%%%%%%%%%%%%%
         %lines
     \draw [line width=1,green](v0) -- (v1);
     \draw [line width=1,green](v0) -- (v2);\draw [line width=1,green](v1) -- (v3);
     \draw [line width=1,green](v3) -- (v5);\draw [line width=1,green](v2) -- (v4);
     \draw [line width=1,green](v5) -- (v6);\draw [line width=1,green](v4) -- (v7);
     \draw [line width=1,green](v6) -- (v10);\draw [line width=1,green](v7) -- (v8);
     \draw [line width=1,green](v8) -- (v9);\draw [line width=1,green](v9) -- (v11);
     \draw [line width=1,green](v10) -- (v15);\draw [line width=1,green](v11) -- (v12);
     \draw [line width=1,green](v12) -- (v13);\draw [line width=1,green](v15) -- (v16);
     \draw [line width=1,green](v13) -- (v14);\draw [line width=1,green](v14) -- (v16);
     %%%%%%%%%%%%%%%%%%%%%%%%%%%%%%%%%%%%%%%%%%%%%%%%%%%%%%%%%%%%%%%%%%%%%%%%%%%%%%%%%%%%%%%%%%%%%%%%%%%%%%%%%%%%%%%
     \draw [dashed,blue](v1) -- (v16);
     \draw [line width=1,blue](v1) -- (v9);\draw [line width=1,blue](v9) -- (v5);
     \draw [line width=1,blue](v5) -- (v12);\draw [line width=1,blue](v16) -- (v2);
     \draw [line width=1,blue](v12) -- (v10);\draw [line width=1,blue](v2) -- (v14);
     \draw [line width=1,blue](v14) -- (v7);\draw [line width=1,blue](v10) -- (v7);
     %%%%%
     \draw [dashed,blue](v4) -- (v13);
     \draw [line width=1,blue](v4) -- (v8);\draw [line width=1,blue](v8) -- (v0);
     \draw [line width=1,blue](v13) -- (v3);\draw [line width=1,blue](v0) -- (v11);
     \draw [line width=1,blue](v3) -- (v15);\draw [line width=1,blue](v15) -- (v6);
     \draw [line width=1,blue](v6) -- (v11);
     %%%%%%%%%%%%%%%%%
      \draw [line width=1,blue](v1) -- (v13);\draw [line width=1,blue](v4) -- (v16);
      %%%%%%%%%%%%%%%%%%%%%%%%%%%%%%%%%%%%%%%%%%%%%%%%%%%%%%%%%%%%%%%%%%%%%%%%%%%%%%%%%%%%%%%%%%%%%%%%%%%%%%%
       \draw [dashed,red](v9) -- (v16);\draw [line width=1,red](v16) -- (v12);
      \draw [dashed,red](v9) -- (v14);\draw [line width=1,red](v12) -- (v14);
      %%%%%%%%%%%%%%%%
      \draw [dashed,red](v8) -- (v13);\draw [line width=1,red](v8) -- (v15);
      \draw [line width=1,red](v13) -- (v11);\draw [line width=1,red](v11) -- (v15);
      %%%%%%%%%%%%%%%%%
      \draw [dashed,red](v1) -- (v10);\draw [line width=1,red](v10) -- (v5);
      \draw [line width=1,red](v5) -- (v7);\draw [line width=1,red](v1) -- (v2);
      \draw [line width=1,red](v7) -- (v2);
      %%%%%%%%%%%%%%%%%
      \draw [dashed,red](v4) -- (v0);\draw [dashed,red,red](v4) -- (v3);
      \draw [line width=1,red](v3) -- (v6);\draw [line width=1,red](v0) -- (v6);
      %%%%%%%%%%%%%%%%%%%%%
      \draw [line width=1,red](v9) -- (v13);\draw [line width=1,red](v8) -- (v16);
      \draw [line width=1,red](v0) -- (v10);\draw [line width=1,red](v1) -- (v4);
      %%%%%%%%%%%%%%%%%%%%%
      \draw [line width=1,red](v9) -- (v3);\draw [line width=1,red](v4) -- (v14);

  %    \node at (v1) {\textbullet};
     \filldraw[black] (v0) circle(2pt);\filldraw[black] (v1) circle(2pt);\filldraw[black] (v2) circle(2pt);\filldraw[black] (v3) circle(2pt);
     \filldraw[black] (v4) circle(2pt);\filldraw[black] (v5) circle(2pt);\filldraw[black] (v6) circle(2pt);\filldraw[black] (v7) circle(2pt);
     \filldraw[black] (v8) circle(2pt);\filldraw[black] (v9) circle(2pt);\filldraw[black] (v10) circle(2pt);\filldraw[black] (v11) circle(2pt);
     \filldraw[black] (v12) circle(2pt);\filldraw[black] (v13) circle(2pt);\filldraw[black] (v14) circle(2pt);\filldraw[black] (v15) circle(2pt);
     \filldraw[black] (v16) circle(2pt);
\end{tikzpicture}
}
\caption {\small A set of  $17\/$ points in general position in the plane with three 1-PHCs where $C_1\/$ in green color, $C_2\/$ in blue color and $C_3\/$ in red color, and the dashed color segments refer to removal joining edge for each cycle}.  \label{three}
\end{figure}
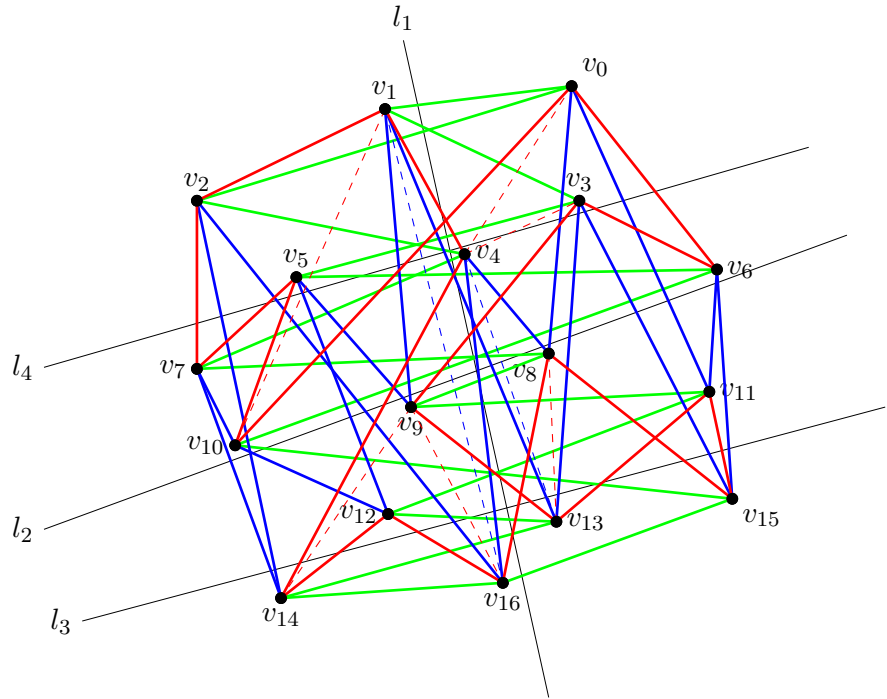

\vspace{1mm} Theorem \ref{the3} is depicted in Figure  \ref{three}.

\begin{thebibliography}{99}

\bibitem{aam:refer} {\sc O. Aichholzer, A. Asinowski, and T. Miltzow}, Disjoint compatibility graph of non-crossing matchings
of points in convex position. {\em Electr. J. Comb.\/}, {\bf 22} (1) (2015) P1.65.

\bibitem{abd:refer} {\sc O. Aichholzer, S. Bereg, A. Dumitrescu, A. G. Olaverri, C. Huemer, F. Hurtado, M. Kano,
A. M´arquez, D. Rappaport, S. Smorodinsky, D. L. Souvaine, J. Urrutia, and D. R. Wood}, Compatible
geometric matchings. {\em Comput. Geom.\/}, {\bf 42}(6-7)(2009) 617 -– 626.

\bibitem{aic:refer} {\sc O. Aichholzer, S. Cabello, R. Fabila-Monroy, D. Flores-Pe\~{n}aloza, T. Hackl, C. Huemer, F. Hurtado and D.R. Wood}, Edge-removal and non-crossing configurations and geometric graphs. {\em Discrete Maths. $\&$ Theorectical Comput. Sci.\/}, {\bf 12} (2010) 75 -- 86.

\bibitem{aki:refer} {\sc S. Akl}, A lower bound on the maximum number of crossing-free Hamilton cycles in a rectilinear drawing of $K_n\/$. {\em Ars Combinatoria \/}, {\bf 7} (1979) 7 -– 18.

\bibitem{bbms:refer} {\sc A. Biniaz, P. Bose, A. Maheshwari and M. Smid}, Packing plane perfect matchings into a point set. {\em Discrete Maths. $\&$ Theorectical Comput. Sci.,\/} {\bf 17} (2015) 119--142.

\bibitem{bk:refer} {\sc F. Bernhart and P. C. Kainen}, The book thickness of a graph. {\em Journal of Combinatorial Theory, Series B\/},
{\bf 27}(3) (1979) 320 -- 331.

\bibitem{bhrw:refer} {\sc P. Bose, F. Hurtado, E. Rivera-Campo, and D. R. Wood}, Partitions of complete geometric graphs into plane trees. {\em  Comput. Geom.\/}, {\bf 34}(2) (2006) 116 -- 125.

\bibitem{dt:refer} {\sc D. Dor and M. Tarsi}, Graph decomposition is NP-complete: A complete proof of Holyer’s
conjecture. {\em SIAM Journal on Computing\/}, {\bf 26}(4) (1997) 1166 –- 1187.

\bibitem{cgh:refer} {\sc M. Claverol, D. Garijo, F. Hurtado, D. Lara and C. Seara}, The alternating path problem revisited. in:{\em Proc. XV Spanish Meeting on Comp. Geom.\/} (2013) 115 -- 118.


\bibitem{cggst:refer} {\sc M. Claverol, A. Garc\'{\i}a, D. Garijo, C. Seara and J. Tejel}, On Hamiltonian alternating cycles and paths. {\em Comput. Geom.} In Press (2017).

\bibitem{dsst:refer} {\sc A. Dumitrescu, A. Schulz, A. Sheffer, and Cs.D. T\'{o}th}, Bounds on the maximum multiplicity
of some common geometric graphs. {\em Proc. 28th Symp. Theoret. Aspects Comput. Sci.\/}, (2011) 637 –- 648.

\bibitem{h:refer} {\sc H. Edelsbrunner}, Algorithms In Combinatorial Geometry. {\em Springer-Verlag, Berlin\/}, (1987).

\bibitem{gjt:refer} {\sc M.R. Garey, D.S. Johnson, and R.E. Tarjan}, The planar Hamiltonian circuit problem is NP-complete. {\em SIAM J. Comput.\/}, {\bf 5} (1976) 704 –- 714.

\bibitem{kky:refer} {\sc A. Kaneko, M. Kano and Y. Yoshimoto}, Alternating Hamiltonian cycles with minimum number of crossings in the plane. {\em Int. J. Comput. Geom. Appl.\/}, {\bf 10} (2000) 73 -– 78.

\bibitem{lot:refer} {\sc M. Li\'{s}kiewicz, M. Ogihara, and S. Toda}, The complexity of counting self-avoiding walks in subgraphs of two-dimensional grids and hypercubes. {\em Theoret. Comput. Sci.\/}, {\bf 304} (2003), 129 –- 156.

\bibitem{cmw:refer} {\sc  C. Lo, J. Matousek, and W. Steiger}, Algorithms for ham-sandwich cuts. {\em Discrete Comput. Geom.\/}, {\bf 11}(4) (1994) 433 –- 452.
\bibitem{m:refer} {\sc T. Motzkin}, Relations between hypersurface cross ratios, and a combinatorial formula for partitions of a polygon, for permanent preponderance, and for non-associative products. {\em Bull. Amer. Math. Soc.\/}, {\bf54}(4) (1948) 352 –- 360.

\bibitem{nm:refer} {\sc M. M. Newborn and W. O. J. Moser}, Optimal crossing-free Hamiltonian circuit drawings of $K_n\/$. {\em J. Comb. Theory, Ser. B\/}, {\bf 29}(1) (1980) 13 -- 26.
    	
\bibitem{ssw:refer} {\sc M. Sharir , A. Sheffe and E. Welzl}, Counting Plane Graphs: Perfect Matchings, Spanning Cycles, and
Kasteleyn’s Technique. {\em J. Combin. Theor. Ser. A\/},  {\bf 120}(4) (2013) 777 -- 794.

\bibitem{sw:refer} {\sc M. Sharir and E. Welzl}, On the number of crossing-free matchings, cycles, and partitions. {\em SIAM J. Comput.\/}, {\bf 36}(3) (2006) 695 -- 720.

\end{thebibliography}
\end{document}